\input ./style/arxiv-vmsta.cfg
\documentclass[numbers,compress,v1.0.1]{vmsta}
\usepackage{vtexbibtags}

\volume{6}
\issue{1}
\pubyear{2019}
\firstpage{3}
\lastpage{12}
\aid{VMSTA127}
\doi{10.15559/18-VMSTA127}
\articletype{research-article}

\startlocaldefs
\allowdisplaybreaks

\endlocaldefs

\begin{document}

\begin{frontmatter}
\pretitle{Research Article}

\title{Probability distributions for the run-and-tumble models with
variable speed and tumbling rate}

\author[a,b]{\inits{L.}\fnms{Luca}~\snm{Angelani}\ead[label=e1]{luca.angelani@roma1.infn.it}}
\author[c]{\inits{R.}\fnms{Roberto}~\snm{Garra}\thanksref{cor1}\ead
[label=e2]{roberto.garra@uniroma1.it}}
\thankstext[type=corresp,id=cor1]{Corresponding author.}

\address[a]{ISC-CNR, \institution{Institute for Complex Systems}, P.le
A. Moro 2, 00185 Rome, \cny{Italy}}
\address[b]{Dipartimento di Fisica, \institution{Sapienza Universit\`a
di Roma}, P.le A. Moro 2, 00185~Rome,~\cny{Italy}}
\address[c]{Dipartimento di Scienze Statistiche, \institution{Sapienza
Universit\`a di Roma}, P.le A. Moro 2, 00185 Rome, \cny{Italy}}


\markboth{L. Angelani, R. Garra}{Run-and-tumble models with variable
speed and tumbling rate}

\begin{abstract}
In this paper we consider a telegraph equation with time-dependent
coefficients, governing the persistent random walk of a particle moving on
the line with a time-varying velocity $c(t)$ and changing direction at
instants distributed according to a non-stationary Poisson distribution
with rate $\lambda(t)$. We show that, under suitable assumptions, we are
able to find the exact form of the probability distribution. We also
consider the space-fractional counterpart of this model, finding the
characteristic function of the related process. A conclusive discussion is
devoted to the potential applications to run-and-tumble models.
\end{abstract}

\received{\sday{16} \smonth{7} \syear{2018}}
\revised{\sday{8} \smonth{11} \syear{2018}}
\accepted{\sday{7} \smonth{12} \syear{2018}}
\publishedonline{\sday{21} \smonth{12} \syear{2018}}
\end{frontmatter}

\section{Introduction}

Many motile bacteria,\index{motile bacteria} such as the common {\it E.coli}, explore the environment
performing run-and-tumble motion \cite{Berg}.
Helicoidal filaments,\index{Helicoidal filaments} called flagella, powered by internal motors
allow the cell to wander around: when flagella rotate counterclockwise
(as seen from behind) the cell performs a straight line motion ({\it run}),
while a clockwise flagellar rotation induces a random reorientation of
the cell body\index{cell body} ({\it tumble}).
In the absence of external force fields or chemicals in the bacterial solution,\index{bacterial solution}
the swim speed $c$ and the rate $\lambda$ at which swimmers change direction
are assumed to be constant in time and space.
In the idealized one-dimensional case the corresponding run-and-tumble
equation of motion
reduces to the usual telegrapher's equation \cite
{Ang_2015,Ang_2017,Ang_2014,mart,sch_1993,were}
%
\begin{equation}
\frac{\partial^2 p}{\partial t^2} + 2 \lambda\frac{\partial p}{\partial t} = c^2
\frac{\partial^2 p}{\partial x^2} .
\end{equation}
%
However, in many interesting real situations swimmers' speed and
tumbling rate\index{tumbling rate} can be
spatial or time dependent quantities.
Recent investigations have demonstrated that the speed of
genetically engineered bacteria, expressing proteorhodopsin protein,
can be tuned by modulating the intensity of an external light field
\cite{Arl_2018,Frangi2018,Tip_2013,Wiz_2017,Wal_2007}.
In such a case one can have a direct control on the swimmers speed by simply
applying a suitable external field.
In particular, time-dependent external fields give rise to time
variable swimmers speed.
Recent investigations have also shown that, for some marine bacteria,\index{marine bacteria}
there is a correlation between the speed and the reorientation frequency.\index{reorientation frequency}
More specifically one observe a linear relationship between the two quantities
in the low-speed regime \cite{Son_2016}.
In such a case it is then appropriate to make the assumption of
proportionality between
$\lambda$ and $c$.

Motivated by these interesting problems, in Section~\ref{sec2}, we provide some
general results
regarding the telegraph equation\index{telegraph equation} with time-dependent parameters $c(t)$
and $\lambda(t)$.
We then analyze the interesting case of proportionality between
velocity\index{velocity} and tumbling rate,\index{tumbling rate} reporting
exact expressions for the probability distribution\index{probability distribution} and the mean square
displacement and discussing the long-time diffusive
behavior for different choice of $c(t)$.

In Section~\ref{sec3} we generalize the above results to the case of the
space-fractional telegraph equation with time-dependent velocity and rate.
Indeed, in the recent literature space and time-fractional
generalizations of the telegraph equations\index{telegraph equation} have attracted the interest
of different authors, see for example \cite
{compte,mirko,fabrizio,giusti,ob,oz}. In \cite{mirko} the relationship
between space-time fractional telegraph equations and time-changed
processes has been discussed. In the recent paper \cite{ma1}, Masoliver
has introduced a fractional persistent random walk,\index{persistent random walk} whose probability
law\index{probability law} is governed by the space-time fractional telegraph equation. The
physical motivation for this kind of generalization is strictly related
to the analysis of sub- and super-diffusive processes, as well as the
telegraph process\index{telegraph process} leads to a ballistic process for short times (and a
classical diffusive one for long times).
We analyze here the space-fractional counterpart of the generalized
telegraph equation studied in Section~\ref{sec2}, finding the characteristic
function\index{characteristic function} of the non-homogeneous fractional telegraph process with
varying velocity.

In a final section we interpret the obtained results in the context of
run-and-tumble models with time-variable swimmers' speed.
In particular, we consider genetically engineered {\it E.coli} bacteria\index{coli bacteria}
whose dynamics is described by run-and-tumble models in which the
value of the speed is controlled by an external field.
We derive the equation of motion in some simple situations, such as
the case of a sudden switch of external fields.

\section{Non-homogeneous telegraph process with time-varying parameters}\label{sec2}

The telegraph process\index{telegraph process} has attracted the interest of many researchers,
starting from the seminal works of Goldstein \cite{golds} and Kac
\cite
{kac}, being a relevant prototype of finite velocity random motion,\index{velocity random motion}
whose probability law\index{probability law} coincides with the fundamental solution of the
telegraph equation.\index{telegraph equation} There is a wide literature about the applications
and generalizations of the telegraph process,\index{telegraph process} we refer to the recent
monograph \cite{rata} for a complete review about this topic. We also
observe that the telegraph equation,\index{telegraph equation} whose origin comes back to the
classical equations of electromagnetism, has been also suggested by
Davydov, Cattaneo and Vernotte as an alternative
to the classical heat equation for diffusion processes with finite
velocity\index{velocity} of propagation, overcoming the so-called paradox of the
infinite velocity\index{infinite velocity} of heat propagation (we refer to the classical review
\cite{Preziosi} and \cite{giusti} about this topic).

A persistent random walk\index{persistent random walk} with a
variable velocity\index{velocity} is studied in \cite{masoliver}, leading to a generalization of the telegraph process.\index{telegraph process}
As discussed in \cite{masoliver} and \cite{were}, in few special cases
the explicit probability law\index{probability law} of this generalized telegraph process\index{telegraph process} can
be found.
Some results about telegraph process\index{telegraph process} with space-varying velocity have
been found in \cite{noi}.

On the other hand, some recent studies have been devoted to a
non-homogeneous version of the telegraph process,\index{telegraph process} where the particle
changes directions at times distributed according to a non-stationary
Poisson distribution with rate $\lambda(t)$. An interesting case was
considered by Iacus in \cite{iacus} and more recently a special case
related to the Euler--Poisson--Darboux equation has been considered in
\cite{eg}, see also \cite{ao}. Moreover, in \cite{claudio}, large
deviations principles have been applied to the non-homogeneous
telegraph process.
More general and relevant models of finite velocity\index{velocity} diffusion
processes are the so-called
L\'evy walks, we refer for example to the recent review \cite{zab}
about this topic.

Here we consider the persistent random walk\index{persistent random walk} of a particle moving on
the line and switching from the time-varying
velocity $c(t)$ to $-c(t)$ at times distributed according to a
non-stationary Poisson distribution with rate $\lambda(t)$. Therefore,
here we consider both the generalizations recently suggested in the
literature and we show that, in a special case, this can help to find
the explicit probability law.\index{probability law}
We assume that $c(t)\in L^1[0,t]$.
According to the classical treatment of the two-direction persistent
random walk\index{random ! walk} given for example by Goldstein \cite{golds} (see also
\cite
{masoliver}), for the description of the position $X(t)$ of the
particle at time $t>0$, we use the probabilities
%
\begin{align}
& a(x,t)dx =P\bigl\{X(t)\in dx, V(t)= c(t)\bigr\},
\\
& b(x,t)dx = P\bigl\{X(t)\in dx, V(t)=- c(t)\bigr\},
\end{align}
satisfying the system of partial differential equations
%
\begin{equation}
\label{sy} %
\begin{cases}
&\displaystyle\frac{\partial a}{\partial t} = -c(t)\frac{\partial
a}{\partial x}+\lambda(t)(b(x,t)-a(x,t)),\\[6pt]
&\displaystyle\frac{\partial b}{\partial t} = c(t)\frac{\partial
b}{\partial x}+\lambda(t)(a(x,t)-b(x,t)),
\end{cases} %
\end{equation}
subject to the initial conditions $a(x,0)= b(x,0)= \frac{1}{2}\delta(x-x_0)$.
The functions $a(x,t)$ and $b(x,t)$ denote the probability density
functions for the position of the random walker\index{random ! walker}
at time $t>0$ while moving respectively in the positive or negative
$x$ direction.
These equations can be simply combined in a single equation for the
total probability $p(x,t) =a(x,t)+b(x,t)$. Let us introduce the
auxiliary function $w(x,t) = a(x,t)-b(x,t)$. By adding and subtracting
the equations in \eqref{sy}, we obtain
%
\begin{equation}
\label{sy1} %
\begin{cases}
&\displaystyle\frac{\partial p}{\partial t} = -c(t)\xch{\frac{\partial
w}{\partial x},}{\frac{\partial
w}{\partial x}}\\[6pt]
&\displaystyle\frac{\partial w}{\partial t} = -c(t)\frac{\partial
p}{\partial x}-2\lambda(t)\xch{w,}{w.}
\end{cases} %
\end{equation}
and finally the following telegraph equation\index{telegraph equation} with time-varying coefficients
%
\begin{equation}
\label{vt} \frac{1}{c(t)}\frac{\partial}{\partial t} \frac{1}{c(t)}
\frac
{\partial
p}{\partial t}+\frac{2\lambda(t)}{c^2(t)}\frac{\partial p}{\partial t} = \frac{\partial^2 p}{\partial x^2}.
\end{equation}
We observe that, from the physical point of view, in the context of the
hyperbolic formulation of the heat wave propagation, equations \eqref
{sy1} are formally equivalent to the heat balance equation with a
time-dependent diffusivity coefficient coupled with a Cattaneo law with
time-varying relaxation.

As pointed out by Masoliver and Weiss in \cite{masoliver}, equations
like \eqref{vt} are generally difficult to be handled \xch{analytically}{analitically}.
However, we observe that, taking $c(t)=c_0 w(t)$, by means of the
change of variable (see also \cite{were})
%
\begin{equation}
\label{change0} \tau= \int_0^t w(s)ds,
\end{equation}
equation \eqref{vt} is reduced to a simpler telegraph-type equation
%
\begin{equation}\label{eq8}
\biggl[\frac{\partial^2}{\partial\tau^2} +2\lambda_{\mathrm{eff}}(\tau )\frac
{\partial}{\partial\tau}
\biggr] p(x, \tau) = c_0^2 \frac{\partial^2
p}{\partial x^2},
\end{equation}
where
%
\begin{equation}
\lambda_{\mathrm{eff}}(\tau) = \frac{\lambda(t(\tau))}{w(t(\tau))}.
\end{equation}
This is a general scheme that allows to find, in some cases, the
explicit form of the probability law\index{probability law} (see also the discussion in
\cite{were}).

We now consider in detail the case $\lambda_{\mathrm{eff}} = \mathrm{const}.$ (i.e.
$\lambda(t)\sim\lambda_0 \ w(t)$) admitting an exact solution. This
means that the rate of changes of directions follows the
velocity-dependence in time. In this case we have that equation \eqref
{vt}, by means of the change of variable \eqref{change0}, is reduced to
%
\begin{equation}\label{eq10}
\biggl[\frac{\partial^2}{\partial\tau^2} +2\lambda_0\frac
{\partial
}{\partial\tau} \biggr]
p(x, \tau) = c_0^2 \frac{\partial^2
p}{\partial x^2},
\end{equation}
corresponding to the classical telegraph equation with velocity\index{velocity} $c_0$
and changing direction rate $\lambda_0$.

Therefore, considering the initial conditions $p(x,0)=\delta(x)$ and
$\partial_\tau p(x,\tau) |_{\tau=0}=0$ and going back to the
variable $t$, we have that the absolutely continuous component of the
probability distribution\index{probability distribution} of the non-homogeneous telegraph process, in
this case is given by
%
\begin{align}\label{eq11}
&P\bigl\{X(t)\in dx\bigr\} = dx \frac{\displaystyle e^{-\lambda_0 \int_0^t
w(s)ds}}{2} \Biggl[\frac{\lambda_0}{c_0}
I_0 \Biggl( \frac{\lambda_0}{c_0} \sqrt{ \Biggl(c_0
\int_0^t w(s) ds \Biggr)^2-x^2}
\Biggr)\nonumber
\\
&+\frac{1}{c_0w(t)}\frac{\partial}{\partial t}I_0 \Biggl(
\frac
{\lambda_0}{c_0} \sqrt{ \Biggl(c_0\int
_0^t w(s) ds \Biggr)^2-x^2}
\Biggr) \Biggr] , \quad|x|< \Biggl(c_0\int_0^t
w(s)ds \Biggr),
\end{align}
where
%
\begin{equation}
I_0(t) = \sum_{k=0}^\infty
\biggl(\frac{t}{2} \biggr)^{2k}\frac{1}{k!^2},
\end{equation}
is a modified Bessel function.
The component of the unconditional probability distribution\index{unconditional probability distribution} that
pertains to the probability of no-changes of direction according to the
Poisson distribution with time-dependent rate $\lambda(t)$ is
concentrated on the boundary $x = \pm\int_0^t c(s)ds$ and it is given by
%
\begin{equation}
P\Biggl\{X(t)= \pm\int_0^t c(s)ds\Biggr\} =
\xch{\frac{e^{-\lambda_0 \int_0^t w(s)ds}}{2}.}{\frac{e^{-\lambda_0 \int_0^t w(s)ds}}{2},}
\end{equation}
Observe that, in the case $c(t) = c_0$ (i.e. $w(t)=1$), we recover the
probability distribution\index{probability distribution} of the classical telegraph process
with rate $\lambda_0$.

An interesting quantity describing the spatial extent of the random
motion\index{random ! motion} is the mean square displacement $r^2$,
i.e. the second moment of the probability distribution\index{probability distribution} (\ref{eq11}):
%
\begin{equation}
r^2(t)= \frac{c_0^2}{2\lambda_0^2} \Biggl[ 2 \lambda_0 \int
_0^t w(s)ds - 1 + e^{- 2 \lambda_0 \int_0^t w(s)ds} \xch{\Biggr].}{\Biggr]}
\end{equation}
The asymptotic behavior of $r^2$, which is linear in $t$ in the
classical persistent random walk,\index{persistent random walk} now \xch{depends}{dependes}
on the long time \xch{behavior}{beahvior} of the velocity\index{velocity} function $c(t)=c_0 w(t)$. We
can distinguish \xch{different}{differrent} regimes.
Assuming a power-law behavior of $w(t)$ at long time, $w(t)\sim
t^{-\beta}$, we have the following cases:
\begin{itemize}
\item[--] $\beta>1$

The integral of $w$, i.e. $\tau(t)$, is finite for $t \to\infty$,
resulting in a
finite asymptotic mean square displacement, $r^2(t) \to \mathrm{const}.$\xch{ The}{. The}
motion is confined in a finite space domain
whose boundaries are at $\xch{x_{B}}{x_{_B}}=\pm\int_0^\infty c_0 w(t) dt$. A
finite stationary probability distribution,\index{stationary probability distribution} given by (\ref{eq11})
in the limit $t \to\infty$, exists.
\item[--] $\beta=1$

In such a case we have logarithmic diffusion, $r^2(t) \sim\ln(t)$
\item[--] $0<\beta<1$

The mean square displacement grows as a power of time,
$r^2(t) \sim t^{\alpha}$,
with $\alpha=1-\beta<1$. The random walk\index{random ! walk} exhibits anomalous diffusion
(subdiffusion).
\item[--] $\beta=0$

In this case the asymptotic velocity is finite, resulting in a linear
time dependence of
the mean square displacement, $r^2(t)\sim t$ (normal diffusion).
\item[--] $\beta<0$

The velocity\index{velocity} grows with time and the random walk\index{random ! walk} is superdiffusive, i.e.
$r^2(t) \sim t^\alpha$, with $\alpha=1-\beta>1$.
\end{itemize}

We can also observe that, assuming an exponential decay $w(t)\sim
e^{-\gamma t}$ a finite stationary probability distribution\index{stationary probability distribution}
exists for $t\rightarrow+\infty$, while if $w(t)$ is a bounded
function, the normal diffusive behaviour is recovered.

\section{The space-fractional telegraph equation with time-varying coefficients}\label{sec3}

The space-fractional telegraph equation was firstly considered by
Orsingher and Zhao in \cite{oz} and more recently studied by Masoliver
in \cite{ma1} in the context of the fractional generalization of the
persistent random walk.\index{persistent random walk}
The relationship between space-time fractional telegraph equations and
time-changed processes have been obtained by D'Ovidio et al. \cite{mirko}.
We here consider the space-fractional telegraph equation with
time-dependent rate and velocity\index{velocity}
%
\begin{equation}
\label{fr1} \frac{1}{c(t)}\frac{\partial}{\partial t} \frac{1}{c(t)}
\frac
{\partial
p}{\partial t}+\frac{2\lambda(t)}{c^2(t)}\frac{\partial p}{\partial t} = \frac{\partial^{2\alpha} p}{\partial|x|^{2\alpha}}, \quad
0<\alpha \leq1.
\end{equation}
The space-fractional derivative appearing in \eqref{fr1} is the Riesz
derivative \cite{kilbas}
%
\begin{equation}
\frac{\partial^{2\alpha} f}{\partial|x|^{2\alpha}}= -\frac
{1}{2\cos
\alpha\pi}\frac{1}{\varGamma(1-2\alpha)}\frac{d}{dx}\int
_{-\infty
}^{+\infty}\frac{f(z)}{|x-z|^{2\alpha}}dz, \quad\alpha\in(0,1),
\end{equation}
whose Fourier transform is given by (see e.g. \cite{mirko} for details)
%
\begin{equation}
\mathcal{F} \biggl[\frac{\partial^{2\alpha} f}{\partial|x|^{2\alpha
}} \biggr](k)= -|k|^{2\alpha}\hat{f}(k),
\end{equation}
where we denote by $\hat{f}(k)$ the Fourier transform of the function $f(x)$.
We here consider in detail the space-fractional counterpart of the
case considered in the previous section, i.e. by taking $c(t) = c_0
w(t)$ and the change of variable $\tau= \int_0^t w(s)ds$, we obtain
the following equation
%
\begin{equation}
\label{fr3} \frac{\partial^2 p}{\partial\tau^2}+2\lambda_{\mathrm{eff}}(\tau)\frac
{\partial p}{\partial\tau}
= c_0^2\frac{\partial^{2\alpha}
p}{\partial
|x|^{2\alpha}},
\end{equation}
where $\lambda_{\mathrm{eff}}(\tau) = \lambda(t(\tau))/w(t(\tau))$.

Considering the special case \mbox{$\lambda(t)\sim\lambda_0 \ w(t)$}, we can
obtain the characteristic function\index{characteristic function} of the non-homogeneous
space-fractional telegraph process with time-var\-ying velocity.
Indeed, we obtain in the Fourier space
%
\begin{equation}
\label{fr2} \frac{\partial^2 \hat{p}}{\partial\tau^2} +2\lambda_0\frac
{\partial
\hat{p}}{\partial\tau} =
-c_0^2|k|^{2\alpha}\hat{p}, \quad 0<\alpha \leq1.
\end{equation}
Therefore, we obtain the characteristic function\index{characteristic function} of the
space-fractional telegraph\break process by means of simple calculations and
going back to the original time variable,
%
\begin{align}
\nonumber
\widehat{p}(k,t) = \frac{e^{-\lambda_0\int_0^t
w(s)ds}}{2}& \biggl[ \biggl(1+
\frac{\lambda_0}{\sqrt{\lambda
_0^2-c_0^2|k|^{2\alpha}}} \biggr)e^{\sqrt{\lambda
_0^2-c_0^2|k|^{2\alpha
}} (\int_0^t w(s)ds )}
\\
&+ \biggl(1-\frac{\lambda_0}{\sqrt{\lambda_0^2-c_0^2|k|^{2\alpha
}}} \biggr) e^{-\sqrt{\lambda_0^2-c_0^2|k|^{2\alpha}} (\int_0^t
w(s)ds
)} \biggr].\label{cf}
\end{align}
The problem to find the inverse Fourier transform of \eqref{cf} seems
to be solvable only in the case $\alpha= 1$ (that leads to the
probability law\index{probability law} of the classical telegraph process).\looseness=1

We can observe that the main features of the space-fractional
telegraph process strongly differ from that of the classical telegraph
process, since it has discontinuous sample paths and it does not
preserve the finite velocity\index{velocity} of propagation.
Indeed, as it was shown by Orsingher and Zhao in \cite{oz}, the random
process related to the space-fractional telegraph equation describes
the one-dimensional motion of a particle which moves forward
and backward performing jumps of random amplitude. This is not
surprising, since the \xch{appearance}{appearence} of the fractional Laplacian is related
to non-locality and leads to almost
surely discontinuous paths.
On the other hand, this model is interesting in the context of the
studies about fractional persistent random walk\index{persistent random walk} models, as fully
discussed by Masoliver in \cite{ma1}.

\section{Discussion: applications to run-and-tumble models}\label{sec4}

We now discuss how the obtained general results can be applied in the
context of run-and-tumble models.

Let us first assume that the tumbling rate\index{tumbling rate} is constant, $\lambda_0$.
This is, for example, the case in which a time-dependent and spatially
homogeneous external field
induces a time-dependent speed $c(t)$ without changing tumbling
processes in genetically engineered bacteria.
In terms of the auxiliary variable $\tau$
the equation of motion turns out to be Eq.~(\ref{eq8})
with a time-dependent effective tumbling rate
%
\begin{equation}
\lambda_{\mathrm{eff}}(\tau)=\xch{\frac{\lambda_0}{w(t(\tau))}.}{\frac{\lambda_0}{w(t(\tau))}}
\end{equation}
A simple interesting case can be analyzed by considering a spatially
uniform light field which is
abruptly switched off at $t=0$. One can assume that, due to finite time
response of the internal processes inside the cell body,\index{cell body}
the swimmer speed exponentially relaxes towards zero
%
\begin{equation}
c(t)=c_0 \exp(-\gamma t) ,
\end{equation}
where $\gamma^{-1}$ is the relaxation time \cite{Arl_2018,Tip_2013}.
In this case one has that
%
\begin{equation}
\lambda_{\mathrm{eff}} = \frac{\lambda_0}{1 -\gamma\tau},
\end{equation}
leading to the partial differential equation
%
\begin{equation}
\label{epd} \biggl[\frac{\partial^2}{\partial\tau^2} +\frac{2\lambda
_0}{1-\gamma
\tau}\frac{\partial}{\partial\tau}
\biggr] p(x, \tau) = c_0^2 \frac
{\partial^2 p}{\partial x^2}.
\end{equation}
We observe that similar equations arise in the analysis of random flights
in higher dimension, see for example \cite{al}.

As mentioned in the Introduction, for some bacteria it has been found
that there is a
proportionality \xch{between}{beweteen} the speed and the reorientation frequency.\index{reorientation frequency}
The assumption $\lambda_{\mathrm{eff}}=\mathrm{const}.$, made in the second part of
Section~\ref{sec2} and
leading to the Eq.~(\ref{eq10}),
is then appropriate for these systems and all the results found in
this approximation apply to this case.
It is still an open question to find, for other microorganisms, the relationship
between tumbling rate\index{tumbling rate} and swim speed.
For example, it would be interesting to investigate
such a issue in the case of genetically engineered bacteria, for which
one could control the bacterial speed by varying the external field
and then
measure the corresponding tumbling rate.\index{tumbling rate}



\end{document}